# Скорость сходимости к предельному режиму в бесконечноканальной системе обслуживания с пуассоновским входящим потоком.

## Головастова Э.А.

### Аннотация.

В данной работе получена оценка скорости сходимости регенерирующего процесса $Q(t)$, значение которого в момент времени $t$ равно числу требований в системе $M\backslash G\backslash\infty$ в этот момент, к предельному (стационарному) режиму исходя из явных формул и также с помощью результатов, касающихся общих регенерирующих процессов.

## 1. Введение.

Рассмотрим систему $M\backslash G\backslash\infty$. Входящий поток - пуассоновский с параметром $\lambda$, времена обслуживания приборов - н.о.р. случайные величины, распределенные по закону $B(x)$; $B(x)$ - непрерывная функция; $\overline{B(x)} = 1 - B(x)$; $b = EB(x) = \int_0^\infty \overline{B(x)}dx$. Заметим, что в данной с.м.о. не может быть очереди, т.к. число приборов бесконечно и прибывшее в систему требование сразу начинает обслуживаться. Пусть $Q(t)$ - число требований в системе в момент времени $t$; $P(Q(0) = 0) = 1$. Это регенерирующий процесс: моменты регенерации - моменты поступления требования в свободную систему; период регенерации есть сумма независимых случайных величин: периода занятости и периода свободного состояния системы. Пусть $T_{рег}$ - случайная величина, задающая период регенерации процесса $Q(t)$, $F(t) = P(T_{рег} < t)$ - ее ф.р.; $T_{зан}$ - период занятости, $G(t) = P(T_{зан} < t)$ - ее ф.р.; $T_{св}$ - период свободного состояния системы,

$P(T_{св} < t) = 1 - e^{-\lambda t}, t > 0.$

Пусть $P_к(t) = P(Q(t) = k)$. Известен следующий результат (см.[4] стр. 206):

$$P_k(t) = \frac{(\rho(t))^k}{k!}e^{-\rho(t)}, k \geq 0 \qquad \rho(t) = \lambda \int_0^t \overline{B(x)}dx$$

$$P_k(\infty) = P_k = \frac{\rho^k}{k!}e^{-\rho} \text{ при } \rho = \lambda b < \infty$$

Пусть:

$$\varphi(t) = \sup_k \left| P_k(t) - P_k \right|$$

Далее мы получим оценку скорости сходимости процесса $Q(t)$ к предельному (стационарному) режиму исходя из явных формул и также с помощью результатов, касающихся общих регенерирующих процессов. В основе последнего - работа [11], в которой приводится явное выражение для распределения периода занятости через распределение времени обслуживания. Это позволяет найти асимптотику распределения периода занятости, а следовательно, и периода регенерации, что впоследствии дает оценку скорости сходимости. Сравним эту оценку с оценкой, полученной на основе точных формул.

## 2. Используемые результаты. Функции распределения времени обслуживания.

Мы рассмотрим следующие классы функций распределения:

1. $B(x)$ имеет легкий хвост, т.е. $B(x)$ убывает не медленнее экспоненты $\overline{B(x)} \leq Ce^{-\lambda t}, \lambda > 0$ или, эквивалентно, $B(x)$ удовлетворяет условию Крамера.

(см.[2] стр. 203 )

$$\exists s_0 > 0: \quad \int_0^\infty e^{s_0 x} dB(x) < \infty.$$

2. В противном функция $B(x)$ имеет тяжелый хвост. Среди функций такого типа мы рассмотрим следующие классы и получим некоторые результаты для их представителей.

    (a) $B(x)$ имеет субэкспоненциальный хвост, т.е.

$$\frac{\overline{B}^{*2}(x)}{\overline{B}(x)} \to 2, \quad x \to \infty,$$

где $f^{*n}(x)$ - $n$−кратная свертка функции $f(x)$ с самой собой.

**Свойство. 1.** *Пусть $\{\xi_i\}, i \geq 1$ н.о.р. случайные величины с субэкспоненциальными хвостами. Случайная величина $\nu$ не зависит от $\xi_i$; $\exists \epsilon > 0 : E(1+\epsilon)^\nu < \infty$. Тогда:*

$$P\left(\sum_{i=1}^\nu \xi_i > x\right) \sim E\nu P(\xi_1 > x) \quad x \to \infty.$$

   (b) $\overline{B(x)}$ - правильно (регулярно) меняющаяся функция т.е.

$$\overline{B(x)} \sim x^{-\alpha} \mathscr{L}(x), \quad x \to \infty \quad 0 < \alpha < \infty$$

Здесь $\mathscr{L}(x)$ - медленно меняющаяся на $\infty$ функция, т.е.

$$\forall t > 0: \quad \frac{\mathscr{L}(tx)}{\mathscr{L}(x)} \to 1, \quad x \to \infty$$

**Свойство. 1.** *Пусть $\{\xi_i\}, i \geq 1$ н.о.р. случайные величины с правильно меняющимися хвостами. Случайная величина $\nu$ не зависит от $\xi_i$; $\exists \epsilon > 0$: $E\nu^{max(1,\alpha+\epsilon)} < \infty$. Тогда:*

$$P(\xi_1 + \xi_2 + \cdots + \xi_\nu > t) \sim E\nu \ P(\xi_1 > t) \quad t \to \infty.$$

Подробнее см. [9].

**Свойство. 2.** *Пусть $\eta = \xi_1 + \xi_2$; $\{\xi_i\}_{i=1,2}$ - независимые неотрицательные случайные величины; $\xi_1$ имеет правильно меняющийся хвост, $\xi_2$ - легкий. Тогда:*

$$P(\eta > t) \sim P(\xi_1 > t) \quad t \to \infty.$$

**Свойство. 3.** *Пусть $Z(x)$ - правильно меняющаяся функция с показателем $\alpha$,*

$$Z_p(x) = \int_0^x y^p Z(y) dy \quad Z_p^*(x) = \int_x^\infty y^p Z(y) dy.$$

## 3. Оценки с помощью точных формул. Общая формула.

$$\left|P_k(t) - P_k\right| = \left|\frac{(\rho(t))^k}{k!}e^{-\rho(t)} - \frac{\rho^k}{k!}e^{-\rho}\right| = |a(t)b(t) - a(\infty)b(\infty)| =$$

$$a(t) = (\rho(t))^k, \ b(t) = \frac{e^{-\rho(t)}}{k!}, \qquad a(\infty) = \rho^k, \ b(\infty) = \frac{e^{-\rho}}{k!}$$

$$= \left|a(t)b(t) - a(\infty)b(t) + a(\infty)b(t) - a(\infty)b(\infty)\right| \leq$$

$$\leq b(t)\left|a(t) - a(\infty)\right| + a(\infty)\left|b(t) - b(\infty)\right| = \frac{e^{-\rho(t)}}{k!}\left|(\rho(t))^k - \rho^k\right| + \frac{\rho^k}{k!}\left|e^{-\rho(t)} - e^{-\rho}\right| =$$

$$= \frac{e^{-\rho(t)}}{k!}\left|(\rho(t) - \rho)\left(\rho(t)^{k-1} + \rho(t)^{k-2}\rho + \cdots + \rho(t)\rho^{k-2} + \rho^{k-1}\right)\right| + e^{-\rho(t)}\frac{\rho^k}{k!}\left|1 - e^{-(\rho-\rho(t))}\right| \leq$$

$$\leq \frac{e^{-\rho(t)}}{k!}k\rho^{k-1}\left|\rho(t) - \rho\right| + e^{-\rho(t)}\frac{\rho^k}{k!}\left|1 - e^{-(\rho-\rho(t))}\right|$$

Используем тот факт, что:

$$\left|1 - e^{-\Delta}\right| < \Delta, \quad \Delta \to 0.$$

И также, что:

$$Arg \max_k \left(\frac{\rho^k}{k!}\right) = [\rho], \quad k \in \mathbb{N} \cup \{0\}.$$

Здесь [ . ]-целая часть.

Тогда:

$$\varphi(t) \leq \sup_k \left(\frac{\rho^{k-1}}{(k-1)!} + \frac{\rho^k}{k!}\right)|\rho - \rho(t)|\,exp\{-\rho(t)\} \leq$$

$$\leq 2\,\frac{\rho^{[\rho]}}{[\rho]!}\left(\rho - \rho(t)\right)exp\{-\rho(t)\}.$$

И в итоге:

$$C_\rho = 2\,\frac{\rho^{[\rho]}}{[\rho]!}$$

$$\varphi(t) \leq C_\rho\left(\rho - \rho(t)\right) = C_\rho\,\lambda\int_t^\infty \overline{B(x)}\,dx. \qquad (1)$$

## 4 Оценки с помощью теории регенерирующих процессов.

Пусть есть регенерирующий, стохастический процесс $\{X(t); t \geq 0\}$ на вероятностном пространстве $(\Sigma, \mathscr{F}, P)$ с измеримым пространством состояний - $(U_X, \mathscr{F}(U_X))$; есть фильтрация - $(\mathscr{F}_{\leq t}, t \geq 0)$, $(\mathscr{F}_{\leq t} \in \mathscr{F})$ и $\{X(t)\}$ измерим по отношению к ней.

**Определение. 4.** *Процесс $\{X(t); t \geq 0\}$ - регенерирующий, если для него существует возрастающая последовательность марковских моментов (т.е. случайных величин, не зависящих от будущего рассматриваемого случайного процесса) $\{\theta_n\}_{n \geq 0}$, $(\theta_0 = 0)$ по отношению к фильтрации $\mathscr{F}_{\leq t}$ таких, что последовательность:*

$$\{\varkappa_n\}_{n=1}^\infty = \{X(t - \theta_{n-1}), \theta_n - \theta_{n-1}, t \in [\theta_{n-1}, \theta_n), \}_{n=1}^\infty$$

*состоит их н.о.р. случайных элементов на $(\Sigma, \mathscr{F}, P)$.*

Случайная величина $\theta_i$ - $i$-ый момент регенерации $\{X(t)\}$; $\tau_i = \theta_i - \theta_{i-1}$ - $i$-ый период регенерации. Пусть:
$$F(t) = P(\tau_1 \leq t), \quad \overline{F(t)} = 1 - F(t)$$
$$\mu = \int_0^\infty \overline{F(t)}dt.$$

Здесь мы полагаем, что $F(t)$ абсолютно непрерывна. Пусть:
$$\Psi(A,t) = P\{X(t) \in A; \tau_1 > t\}; \quad A \in \mathscr{F}(U_X) \quad t \geq 0$$
$$P(A,t) = P(X(t) \in A) \quad A \in \mathscr{F}(U_X) \quad t \geq 0$$

Известно, что (см. [10]):
$$\lim_{t \to \infty} P(A,t) = \frac{1}{\mu} \int_0^\infty \Psi(A,y)dy = P(A),$$

**Теорема. 2.** *Пусть функция $F(t)$ имеет тяжелый хвост; $\mu < \infty$;*
$$V(t) = \int_t^\infty \overline{F(y)}dy, \quad f(t,A) = V^{-1}(t) \int_t^\infty \Psi(A,y)dy.$$

*Пусть также функция $F(t)$ такая, что:*
$$\varlimsup_{t < \infty} \frac{V(t)}{V(2t)} < \infty. \tag{5}$$

*Тогда будет верно:*
$$P(A,t) - P(A) = \mu^{-1}V(t)(P(A) - f(t,A)) + O(u(t)),$$
*где* $u(t) = \overline{F(t)} \int_0^t V(y)dy.$

Подробнее см. [7].
Известно, что (см. [6]):
$$|P(A) - f(t,A)| \leq 1.$$

Тогда в случае тяжелых хвостов у $F(t)$ при выполнении условий теоремы 2 и также, если $u(t) = o(V(t))$, будет для $\forall \epsilon > 0 \; \exists t_\epsilon: \; \forall t > t_\epsilon$:
$$|P(A,t) - P(A)| \leq \mu^{-1}(1+\epsilon)V(t), \quad t \to \infty. \tag{6}$$

В качестве процесса $\{X(t)\}$ мы рассмотрим процесс $Q(t)$, описанный в разделе 1. Тогда $A = \{\omega : Q(t,\omega) = k\} \; k \geq 0; \; P(A,t) = P_k(t), P(A) = P_k$. Также:
$$\overline{F(x)} = P(T_{\text{рег}} > x) = P(T_{\text{зан}} + T_{\text{св}} > x) = P(T_{\text{св}} > x) + \int_0^x P(T_{\text{зан}} > x-y)P(T_{\text{св}} \in dy) =$$
$$= e^{-\lambda x} + \lambda \int_0^x \overline{G(x-y)} \, e^{-\lambda y}dy \tag{7}$$

Т.о. из вида 7 можно предположить, что асимптотика хвоста ф.р. периода регенерации в случае тяжелых хвостов у ф.р. периода занятости будет определяться именно ф.р. периода занятости. Пусть:
$$C(t) = 1 - exp\left\{-\lambda \int_0^t \overline{B(x)}dx\right\}, \quad t \geq 0$$
$$c(t) = C'(t) = \lambda \overline{B(t)}(1 - C(t)).$$

**Теорема. 3.**
$$G(t) = 1 - \lambda^{-1} \sum_{n=1}^\infty c^{*n}(t).$$

Подробнее см. [11].

## 4.1 Легкие хвосты.

Найдем оценку хвоста ф.р. периода занятости. Для этого воспользуемся теоремой 3. В нашем случае:
$$c(t) \le \lambda \overline{B(t)} \le \lambda C e^{-\mu t}.$$

Вычислим:
$$\int_0^\infty c(t)\, dt = \int_0^\infty \lambda \overline{B(t)}\, exp\left\{-\lambda \int_0^t \overline{B(y)} dy\right\} dt =$$
$$= \int_0^\infty exp\left\{-\lambda \int_0^t \overline{B(y)} dy\right\} d\left(\int_0^t \overline{B(y)} dy\right) =$$
$$= \int_0^{\lambda b} e^{-z} dz = 1 - e^{-\lambda b}.$$

Введем:
$$v(t) = \frac{c(t)}{1 - e^{-\lambda b}}.$$

Тогда $v(t)$ есть функция плотности некоторой случайной величины. Теперь:
$$\overline{G(t)} = \lambda^{-1} \sum_{n=1}^\infty (1-e^{-\lambda b})^n v^{*n}(t) = \lambda^{-1}(1-e^{-\lambda b}) e^{\lambda b} \sum_{n=1}^\infty e^{-\lambda b}(1-e^{-\lambda b})^{n-1} v^{*n}(t).$$

Пусть $\{\xi_i\}, i \ge 1$ н.о.р. случайные величины с ф.р. $\int_0^t v(x)dx$. Тогда:
$$\int_t^\infty \overline{G(x)} dx = \frac{e^{\lambda b} - 1}{\lambda} P(\xi_1 + \xi_2 + \cdots + \xi_\nu > t),$$

где $\nu$ имеет геометрическое распределение, $P(\nu = j) = e^{-\lambda b}(1-e^{-\lambda b})^{j-1}$, $j = 1, 2, \ldots$; $\nu$ не зависит от $\{\xi_i\}$.

Т.к. $B(x)$ имеет легкий хвост, то из условия Крамера следует, что
$$\exists\, s_* = sup\{s:\ \int_0^\infty e^{sx} d\, B(x) < \infty\} > 0.$$
$$v^*(s) = \int_0^\infty e^{sx} v(x) dx = \int_0^\infty e^{sx} \frac{\lambda \overline{B(x)}}{1 - e^{-\lambda b}}\, exp\left\{-\lambda \int_0^x \overline{B(y)} dy\right\} dx \le$$
$$\le \frac{\lambda}{1 - e^{-\lambda b}} \int_0^\infty e^{sx} \overline{B(x)} dx.$$

Т.к. $\overline{B(x)} \le Ce^{-\mu x}$, то:
$$v^*(s) \le \frac{K}{\mu - s} < \infty \text{ при } \forall\, 0 < s < \mu.$$

Поэтому ф.р. $\xi_i$ удовлетворяет условию Крамера, т.е. она имеет легкий хвост. Тогда
$$\exists\, s_0 = sup\{s:\ v^*(s) < \infty\} > 0.$$

Докажем теперь, что ф.р. $\sum_{i=1}^\nu \xi_i$ тоже удовлетворяет условию Крамера.
$$\int_0^\infty e^{sx} dP(\xi_1 + \xi_2 + \cdots + \xi_\nu < x) = \frac{pv^*(s)}{1 - qv^*(s)}, \quad p = e^{-\lambda b},\ q = 1 - p.$$

Нужно, чтобы последняя дробь была конечна и положительна. Она будет таковой при $v^*(s) < \frac{1}{q}$. Пусть $s_1 = sup\{s : v^*(s) < \frac{1}{q}\} \leq s_0$. Такое $s_1$ существует, т.к. $v^*(s)$ - непрерывная, монотонно возрастающая функция, и тогда $s_1$ - решение уравнения $v^*(s) = \frac{1}{q}$.

Тогда $\forall s < s_1$ интеграл из условия Крамера для ф.р. $\sum_{i=1}^{\nu} \xi_i$ будет иметь нужные свойства. Т.е $\exists \alpha > 0, C_\alpha > 0$:
$$P(\xi_1 + \xi_2 + \cdots + \xi_\nu > t) \leq C_\alpha e^{-\alpha t}.$$

Тогда:
$$\int_t^\infty \overline{G(x)} dx \leq K_1 e^{-\alpha t} \ \forall x \geq 0.$$

Тогда у ф.р. периода занятости есть конечное м.о.:
$$m = \int_0^\infty \overline{G(x)} dx \leq K_1.$$

$\int_0^t \overline{G(x)} dx$ - монотонно возрастающая к $m$ функция, т.к. $\overline{G(x)}$ неотрицательна при $x \geq 0$.

Тогда $\frac{\int_0^t \overline{G(x)} dx}{m}$ - можно считать ф.р. некоторой случайной величины, причем эта ф.р. будет иметь легкий хвост, т.е. для нее выполнено условие Крамера. Докажем, что для $G(x)$ оно тоже выполнено.

$$\int_0^\infty e^{sx} d\frac{\int_0^x \overline{G(y)} dy}{m} = \int_0^\infty e^{sx} \frac{\overline{G(x)}}{m} dx =$$
$$= \frac{1}{s} e^{sx} \frac{\overline{G(x)}}{m} \bigg|_0^\infty - \frac{1}{s} \int_0^\infty e^{sx} d \frac{\overline{G(x)}}{m} =$$
$$= -\frac{1}{sm} + \frac{1}{sm} \int_0^\infty e^{sx} d\ G(x).$$

Т.е. получили, что:
$$\int_0^\infty e^{sx} d\ G(x) = 1 + sm \int_0^\infty e^{sx} d\frac{\int_0^x \overline{G(y)} dy}{m}.$$

Т.е. интеграл из условия Крамера для $G(x)$ конечен для всех $0 < s$, при которых выполнено условие Крамера для $\frac{\int_0^t \overline{G(x)} dx}{m}$. Поэтому $G(x)$ имеет легкий хвост.

Тогда из 7):
$$\overline{F(t)} \leq e^{-\lambda t} + \lambda \int_0^t K_2\ e^{-\alpha(t-y)}\ e^{-\lambda y} dy = e^{-\lambda t} + \frac{K_2 \lambda}{\alpha - \lambda} e^{-\alpha t}\ e^{(\alpha - \lambda) y} \bigg|_0^t = K_4 e^{-\lambda t} + K_3 e^{-\alpha t}.$$

Т.е. $F(t)$ имеет легкий хвост, поэтому можно воспользоваться результатами Theorem 1. из [6].

**Теорема. 4.** *Пусть $F(t)$ абсолютно непрерывная функция и для некоторого $s_0 > 0$*
$$\int_0^\infty e^{sx} dF(x) < \infty.$$

*Тогда:*
$$\exists \delta, \quad 0 < C_\delta < \infty : \ \varphi(t) < C_\delta e^{-\delta t} \quad t \geq 0.$$

Тогда и в нашем случае $\exists \delta, \quad 0 < C_\delta < \infty$:
$$\varphi(t) < C_\delta e^{-\delta t} \quad t \geq 0.$$

Заметим, что в оценке 2 имелся явный вид констант $C_\delta$ и $\delta$, тогда как в последнем мы можем сказать только про их существование.

## 4.2 Тяжелые хвосты. Правильно меняющиеся функции.

Положим здесь $\alpha > 2$.

Далее проведем рассуждения, аналогичные рассуждениям раздела 4.1. Согласно [11],

$$\overline{G(t)} = \lambda^{-1} \sum_{k=1}^{\infty} c^{*k}(t).$$

Определяем функцию:

$$v(t) = \frac{c(t)}{1 - e^{-\lambda b}},$$

которая есть плотность некоторой случайной величины. Тогда:

$$\overline{G(t)} = \lambda^{-1}(1 - e^{-\lambda b})e^{\lambda b} \sum_{k=1}^{\infty} e^{-\lambda b}(1 - e^{-\lambda b})^{k-1} v^{*k}(t).$$

Пусть $\{\xi_i\}, i \geq 1$ н.о.р. случайные величины с ф.р. $\int_0^t v(x)dx$. Тогда:

$$\int_t^{\infty} \overline{G(x)}dx = \frac{e^{\lambda b}(1 - e^{-\lambda b})}{\lambda} P(\xi_1 + \xi_2 + \cdots + \xi_\nu > t),$$

где $\nu$ имеет геометрическое распределение, $P(\nu = j) = e^{-\lambda b}(1 - e^{-\lambda b})^{j-1}, j = 1, 2, \ldots$; $\nu$ не зависит от $\{\xi_i\}$; $E\nu = e^{\lambda b}$.

Т.к.

$$v(t) = \frac{\lambda}{1 - e^{-\lambda b}} \overline{B(t)} \, exp\left\{-\lambda \int_0^t \overline{B(x)}dx\right\} \sim \frac{\lambda e^{-\lambda b}}{1 - e^{-\lambda b}} \, t^{-\alpha} \mathscr{L}(t) \quad t \to \infty,$$

т.е. плотность $\xi_i$ есть правильно меняющаяся функция. Тогда с использованием свойства 3) для правильно меняющихся функций и что $\alpha > 2$:

$$\int_t^{\infty} v(x)dx \sim \frac{tv(t)}{|\alpha + 1|},$$

получаем, что ф.р. $\xi_i$ имеет правильно меняющийся хвост. Поэтому можем воспользоваться свойством 2. Тогда:

$$P(\xi_1 + \xi_2 + \cdots + \xi_\nu > t) \sim E\nu \, P(\xi_1 > t).$$

Поэтому:

$$\int_t^{\infty} \overline{G(x)}dx \sim \frac{e^{2\lambda b}(1 - e^{-\lambda b})}{\lambda} P(\xi_1 > t) = \frac{e^{2\lambda b}(1 - e^{-\lambda b})}{\lambda} \int_t^{\infty} v(x)dx.$$

По правилу Лопиталя:

$$1 = \lim_{t \to \infty} \frac{\int_t^{\infty} \overline{G(x)}dx}{K_1 \int_t^{\infty} v(x)dx} = \lim_{t \to \infty} \frac{\overline{G(t)}}{K_1 v(t)}.$$

Т.е. ф.р. периода регенерации имеет правильно меняющийся хвост:

$$\overline{G(x)} \sim e^{\lambda b} \, \overline{B(x)}.$$

Случайные величины $T_{\text{зан}}$ и $T_{\text{св}}$ удовлетворяют условиям свойства 1, поэтому:

$$\overline{F(t)} \sim \overline{G(t)} \quad t \to \infty.$$

Т.е. хвост ф.р. периода регенерации тоже правильно меняется.

Из теоремы 5 получаем, что здесь у ф.р. периода регенерации есть конечное м.о..

С учетом свойства 3) правильно меняющихся функций имеем:

$$V(t) \sim \int_t^{\infty} K\overline{B(y)}dy \sim \frac{K}{|\alpha + 1|} t\overline{B(t)}.$$

Т.к. у нас $\alpha > 2$, то мы получили, что $V(t)$ тоже правильно меняется.

Т.к. у нас $\alpha > 2$, то мы получили, что $V(t)$ тоже правильно меняется.

$$u(t) \sim K_2 \overline{B(t)} \int_0^t V(y) dy \sim K_3 t^2 \overline{B(t)}^2.$$

Учитывая вышесказанное, можно воспользоваться результатами Theorem 2. из [6]:

**Теорема. 6.** *Пусть*

$$\overline{F(x)} \sim x^{-\alpha} \mathscr{L}(x), \quad x \to \infty,$$

*где $\alpha > 2$ и $\mathscr{L}(x)$ медленно меняющаяся функция. Тогда для любого $\epsilon > 0$ существует $t_\epsilon$ такое, что:*

$$\varphi(t) \leq \frac{(1+\epsilon)\mathscr{L}(t)}{\mu t^{\alpha-1}} \text{ при } t \geq t_\epsilon.$$

Тогда в нашем случае $\forall \epsilon > 0 \; \exists t_\epsilon : \forall t \geq t_\epsilon$

$$\varphi(t) \leq \frac{e^{\lambda b}(1+\epsilon)\mathscr{L}(t)}{\mu(\alpha+1) \; t^{\alpha-1}}, \quad t \to \infty.$$

# Список литературы